\documentclass[11pt,a4paper]{article}

\usepackage[margin=1in]{geometry}
\usepackage{amsmath,amsthm,amssymb}
\usepackage[final]{pdfpages}
\usepackage{hyperref}

\theoremstyle{remark}

\newcommand{\N}{\mathbb{N}}

\newcommand{\R}{\mathbb{R}}
 
\title{On Enhancing the Dissipative Behavior of Active Flux Advection Schemes}

\author{Christian Klingenberg\thanks{Institute of Mathematics, University of W\"urzburg, W\"urzburg, Germany} , Simon Krotsch\footnotemark[1] , Philip Roe\thanks{Formerly: Department of Aerospace Engineering, University of Michigan, Ann Arbor, USA
}}

\date{\vspace{-10mm}}
 
\begin{document}

\maketitle

\begin{abstract}
	In this work, the traditional third-order Active Flux advection scheme is modified by reformulating the method and introducing four additional parameters $R$, $S$, $T$, and $U$ that allow the numerical dissipation and dispersion of the scheme to be adjusted. The influence of these parameters is studied using Fourier analysis and parameter choices that lead to schemes with improved dissipative properties are identified. These improvements are validated by numerical experiments.
\end{abstract}

\section{Introduction}

Some 15 years ago, Phil Roe introduced the Active Flux method \cite{eymann_roe_2011}, \cite{eymann_roe_2011_2}, \cite{eymann_roe_2013}, a numerical method for conservation laws inspired by an idea of Bram van Leer \cite{van_leer_1977}. It is a finite volume method. The degrees of freedom of the traditional third-order Active Flux method are cell averages and point values located at cell interfaces. These point values are evolved in time independently of the averages and are shared between two adjacent cells, resulting in a globally continuous reconstruction and, hence, a Riemann-solver-free method.

On April 26, 2026, Phil Roe passed away. Until the end of his life, he pursued the development of the Active Flux method with unwavering commitment. Examples of his contributions to the development of Active Flux forward include an extension of the method to advection-diffusion problems in  \cite{nishikawa_roe_2016}. An example of the connection to the superconvergence of DG methods was found in \cite{roe_superconvergence_2017}, and a practical comparison between the two approaches can be found in \cite{roe_maeng_fan_2018}. In \cite{barsukow_etal_2026}, this comparison has now been expanded upon. Furthermore, a high-order version was introduced in \cite{samani_roe_2023}. Phil Roe's Active Flux work also inspired further contributions by other groups; some examples are \cite{Abgrall2026}, \cite{helzel_2023},
\cite{barsukow_etal_2019}, 
\cite{barsukow_arxiv_2025},
and \cite{Duan_2025}.

The ICOSAHOM 2025 in Montr\'eal was the last conference Phil Roe attended in person. In memory of Phil Roe, we present one of his last ideas concerning the Active Flux method.
First, the traditional third-order Active Flux method is reviewed and reformulated as an upwind method with correction terms depending on positive parameters. The dissipative and dispersive properties of this parameterized formulation are analyzed, and parameter choices are determined that yield methods with improved dissipative properties. Phil Roe called the method with impressively low dissipation the Super-Duper Method. Finally, numerical experiments are conducted to further investigate the proposed methods and verify their enhanced properties. 

\section{The Active Flux method}

We consider the third-order Active Flux method to solve initial value problems for the one-dimensional linear advection equation
\begin{align*}
	\frac{\partial}{\partial t} q(t,x) + a \frac{\partial}{\partial x} q(t,x) = 0, \quad x \in \mathbb{R}, \ t \in \mathbb{R}_+, 
\end{align*}
where $a > 0$. 
Assume the computational domain is divided into grid cells $\mathcal{C}_i = [x_{i-\frac{1}{2}}, x_{i+\frac{1}{2}}]$ of equal length $\Delta x > 0$. The midpoint of the cell $\mathcal{C}_i$ is denoted by $x_i$.
The degrees of freedom of the traditional third-order Active Flux method are cell averages 
\begin{align*}
	Q_i(t) \approx \frac{1}{\Delta x}\int_{\mathcal{C}_i}q(t, x) \ \mathrm{d}x
\end{align*}
and point values
\begin{align*}
	q_{i+\frac{1}{2}}(t) \approx q(t, x_{i+\frac{1}{2}})
\end{align*}
at the interfaces of the cells. Since point values are shared between adjacent cells, the degrees of freedom associated with the cell $\mathcal{C}_i$ are a point value and a cell average. If these values are given at some time $t^n \geq 0$, $n \in \N_0$, for all cells, the first step of the method is to fit a quadratic polynomial to the three degrees of freedom accessible per cell. An easy way to state this quadratic polynomial, which we call the reconstruction, is to consider a reference cell $\left[-\frac{\Delta x}{2}, \frac{\Delta x}{2}\right]$. The local reconstruction in cell $\mathcal{C}_i$ is defined as {\color{black}$q_{\mathrm{recon},i} : [0, \infty) \times  \left[-\frac{\Delta x}{2}, \frac{\Delta x}{2}\right]  \rightarrow \R$},
\begin{align*}
	q_{\mathrm{recon},i} (t^n, \xi) = & \ \frac{6Q_i(t^n) - q_{i+\frac{1}{2}}(t^n) - q_{i-\frac{1}{2}}(t^n)}{4} \\
    & \ + \frac{q_{i+\frac{1}{2}}(t^n) - q_{i-\frac{1}{2}}(t^n)}{\Delta x}\xi +3 \frac{q_{i-\frac{1}{2}}(t^n) + q_{i+\frac{1}{2}}(t^n) - 2Q_i(t^n)}{\Delta x^2}\xi^2,
\end{align*} 
where $\xi = x - x_i$ for $x \in \mathcal{C}_i$.
One can verify that this function interpolates the degrees of freedom, i.e., 
\begin{align*}
	q_{\mathrm{recon},i}\left(t^n, \pm\frac{\Delta x}{2}\right) = q_{i \pm \frac{1}{2}}(t^n), && \frac{1}{\Delta x}\int_{-\frac{\Delta x}{2}}^{\frac{\Delta x}{2}}q_{\mathrm{recon},i} (t^n, \xi) \ \mathrm{d}\xi = Q_i(t^n).
\end{align*} 
From the integral form of the conservation law, one can derive an update for the cell average: 
\begin{align*}
	Q_i(t^{n+1}) = Q_i(t^{n}) - \frac{a \Delta t}{\Delta x}\bigg(\overline q_{i+\frac{1}{2}}^{n+\frac{1}{2}} - \overline q_{i-\frac{1}{2}}^{n+\frac{1}{2}}\bigg),
\end{align*}
where
\begin{align*}
	 \overline q_{i+\frac{1}{2}}^{n+\frac{1}{2}} = \frac{1}{\Delta t}\int_{t^n}^{t^{n+1}} q_{i+\frac{1}{2}}(t) \ \mathrm{d}t \approx \frac{1}{\Delta t}\int_{t^n}^{t^{n+1}} q(t, x_{i+\frac{1}{2}}) \ \mathrm{d}t
\end{align*}
and $\Delta t = t^{n+1}-t^n$. To evaluate the integral, we use Simpson's rule: 
\begin{equation*}
	\overline q_{i+\frac{1}{2}}^{n+\frac{1}{2}} = \frac{1}{6}\bigg( q_{i+\frac{1}{2}}(t^n) + 4 q_{i+\frac{1}{2}}(t^n + \frac{1}{2}\Delta t) + q_{i+\frac{1}{2}}(t^{n} +  \Delta t)\bigg).
\end{equation*}
The values $q_{i+\frac{1}{2}}(t^{n} + l\Delta t)$ for $l = \frac{1}{2}$ and $l = 1$ are obtained by solving the initial value problem at the locations of the point values with the reconstruction as initial data. For the linear advection equation with $a > 0$, the solution is
\begin{equation*}
	q_{i+\frac{1}{2}}(t^{n} + l\Delta t) =  q_{\mathrm{recon},i}(t^n, \tfrac{\Delta x}{2} - al\Delta t)
\end{equation*}
under the condition
\begin{align*}
	\nu = \frac{a\Delta t}{\Delta x} \leq 1.
\end{align*}

 Since the update of the average is exact, this results in a third-order accurate method. For simplicity, define $Q_i^n = Q_i(t^n)$ and $q_{i \pm \frac{1}{2}}^n =q_{i \pm \frac{1}{2}}(t^n)$. The update formulas for the point values and the cell averages can be written explicitly:
\begin{align*}
	q^{n+1}_{i+\frac{1}{2}} &= (1-\nu ) q_{i+\frac{1}{2}}^n + \nu q_{i-\frac{1}{2}}^n - 3\nu (1-\nu) \left(\left(q_{i+\frac{1}{2}}^n - Q_i^n\right)  - \left(Q_i^n - q_{i-\frac{1}{2}}^n\right)\right), \\
	\overline q^{n+\frac{1}{2}}_{i+\frac{1}{2}} &= Q_i^n + (1-\nu) \left((1-\nu) \left(q_{i+\frac{1}{2}}^n - Q_i^n\right)  + \nu\left (Q_i^n - q_{i-\frac{1}{2}}^n\right)\right), \\
	Q_i^{n+1} &= Q_i^n - \nu \left(\overline q^{n+\frac{1}{2}}_{i+\frac{1}{2}} - \overline q^{n+\frac{1}{2}}_{i-\frac{1}{2}}\right).
\end{align*}
These explicit formulas can be interpreted as a first-order upwind scheme applied separately to the cell average and the point values with additional correction terms. We can write these correction terms in a more general way by introducing parameters $R, S, T , U> 0$:
\begin{align}
	\label{eq:scheme_mod_1}
	q^{n+1}_{i+\frac{1}{2}} &= (1-\nu ) q_{i+\frac{1}{2}}^n + \nu q_{i-\frac{1}{2}}^n - \nu (1-\nu) (R(q_{i+\frac{1}{2}}^n - Q_i^n)  - S(Q_i^n - q_{i-\frac{1}{2}}^n)), \\
	\label{eq:scheme_mod_2}
	\overline q^{n+\frac{1}{2}}_{i+\frac{1}{2}} &= Q_i^n + (1-\nu) (T (q_{i+\frac{1}{2}}^n - Q_i^n)  + U (Q_i^n - q_{i-\frac{1}{2}}^n))
\end{align}
By choosing $R = S = 3$, $T=(1-\nu)$, and $U = \nu$, one recovers the classical third-order scheme. 

\section{Spectral error analysis}

The degrees of freedom of the scheme can be regarded as functions defined on a finite grid. Such functions can be decomposed into finite Fourier series in space at each time step $t^n$:
\begin{equation*}
	Q_i^n  = \sum_{\theta} \hat U ^n (\theta) \mathrm e^{-\frac{1}{2}I\theta} \mathrm e^{I(i+\frac{1}{2})\theta}
\end{equation*}
and
\begin{equation*}
	q_{i+\frac{1}{2}}^n  = \sum_{\theta} \hat q^n(\theta) \mathrm e^{I(i+\frac{1}{2})\theta},
\end{equation*}
where $I$ is the imaginary unit and $\theta = k\Delta x$ the phase angle with wavenumber $k$. For convenience, we define a new amplitude for the cell averages, $\hat Q^n(\theta) = \hat U ^n (\theta) \mathrm e^{-\frac{1}{2}I\theta}$. Assume that the average and point values are given by a single harmonic, i.e., 
\begin{align*}
	Q_i^n  = \hat Q ^n (\theta) \mathrm e^{I(i+\frac{1}{2})\theta}.
\end{align*}
The exact time evolution of the point values and cell averages is given by
\begin{align*}
	q\left(t^{n+1} ,x_{i+\frac{1}{2}}\right) &= q\left(t^{n} ,x_{i+\frac{1}{2}}\right) \mathrm e^{-I\theta\nu}, \\
	\frac{1}{\Delta x} \int_{\mathcal{C}_i}q(t^{n+1},x)\ \mathrm{d}x &= \frac{1}{\Delta x} \int_{\mathcal{C}_i}q(t^{n},x) \ \mathrm{d}x \  \mathrm e^{-I\theta\nu}.
\end{align*} 

Evolving the degrees of freedom in time with the parameterized version of the Active Flux scheme yields
	\begin{align*}
		\hat q^{n+1}   = \ & (1-\nu ) \hat q^{n}   + \nu \hat q^{n}  \mathrm e^{-I\theta} - \nu (1-\nu) \left(R\left(\hat q^{n}  - \hat Q^n\right)  - S\left(\hat Q^n - \hat q^{n}  \mathrm e^{-I\theta}\right)\right) 
	\end{align*}
and 
	\begin{align*}
		\hat Q^{n+1} = & \ \hat Q^n - \nu \bigg(\hat Q^n +(1-\nu)\left(T\left(\hat q^n - \hat Q^n \right) + U\left(\hat Q^n  - \hat q^n \mathrm{e}^{-I\theta}\right)\right) \\
		&\ - \hat Q^n  e^{-I\theta} + (1-\nu)\left(\left(T(\hat q^n - \hat Q^n \right) + U\left(\hat Q^n  - \hat q^n \mathrm{e}^{-I\theta}\right)\right)\mathrm{e}^{-I\theta} \bigg).
	\end{align*}
For simplicity, we omit the dependence on $\theta$. The amplification matrix $A \in \mathbb{C}^{2 \times 2}$ of this scheme, i.e., the matrix that satisfies
    \begin{align*}
	\begin{pmatrix}
			\hat Q^{n+1}  \\
			\hat q^{n+1}
		\end{pmatrix} = A \begin{pmatrix}
			\hat Q^n  \\
			\hat q^n
		\end{pmatrix},
    \end{align*}
is
	\begin{align*}
		A = \begin{pmatrix}
			1 -\nu(1 + (1-\nu)(U-T))\left(1-\mathrm{e}^{-I\theta}\right) & \nu(1-\nu)\left(U\mathrm{e}^{-I\theta} - T\right)\left(1-\mathrm{e}^{-I\theta}\right) \\
			\nu(1-\nu)(R + S) & (1-\nu)(1-\nu R) + \nu(1 -  (1-\nu)S)\mathrm e^{-I\theta}
		\end{pmatrix}.
	\end{align*}
Mathematica \cite{mathematica_2025} allows us to compute the eigenvalues of $A$ as power series in $\theta$. The first two terms in the expansion of the first eigenvalue, which we will also call the principal eigenvalue, are $1$ and $-I \nu \theta$, while the higher-order terms depend on the positive parameters $R$, $S$, $T$, and $U$. In contrast, all terms in the expansion of the second (spurious) eigenvalue depend on these parameters. The principal eigenvalue describes the behavior of low frequencies and the spurious eigenvalue governs the behavior of high frequencies. The exact eigenvalue is
\begin{align*}
	\lambda_{\mathrm{exact}}{\color{black}(\nu, \theta)} = \mathrm{e}^{-I\nu\theta} = \sum_{j=0}^\infty \frac{(-I\nu\theta)^j}{j!}.
\end{align*}
Hence, the first two terms of the principal eigenvalue match the terms of the exact eigenvalue. We can determine the parameters $R,S,T,U$ such that the principal eigenvalue is exact up to the highest possible order of $\theta$. This can be done by employing Mathematica \cite{mathematica_2025}. The resulting conditions on the parameters to obtain a second-order correct principal eigenvalue are:
\begin{align*}
	U =  \frac{R-2 S T+S}{2 R},
\end{align*}
or
\begin{align*}
	R =  0 \quad \text{and} \quad T = \frac{1}{2}.
\end{align*}
To obtain a third-order correct eigenvalue, the parameters must satisfy
\begin{align}
	T =  -\frac{1}{3} (\nu +1) R+\frac{R^2}{R+S}+\frac{1}{2} \quad \text{and} \quad U =  \frac{1}{6}
   \left(2 S \left(\nu -\frac{3 R}{R+S}+1\right)+3\right).
   \label{eq:cond_2}
\end{align}
If conditions \eqref{eq:cond_2} are satisfied, the difference between the exact and principal eigenvalue is
\begin{align*}
	\lambda_1{\color{black}(\nu, \theta)} - \lambda_{\mathrm{exact}}{\color{black}(\nu, \theta)} = \frac{1}{72} (\nu -1) \nu  \left((\nu
   -2) (\nu +1)+\frac{18}{R+S}\right) \theta ^4+O\left(\theta ^5\right).
\end{align*}
If $R = S = 3$, conditions \eqref{eq:cond_2} reduce to  $T = 1 - \nu $ and $U = \nu$, which are the parameters corresponding to the traditional Active Flux method. The conditions for a fourth-order correct eigenvalue are
\begin{align}
\begin{split}
	S &= \frac{18}{-\nu ^2+\nu +2}-R, \\
	T &= \frac{1}{18} (9-(\nu +1) R ((\nu -2) R+6)), \\
	U &= \frac{1}{18} \left(-\left((\nu -2) (\nu +1) R^2\right)-6 (\nu +4)
   R+\frac{9 (\nu -14)}{\nu -2}\right).
 \end{split}
 \label{eq:cond_3}
\end{align}
The error of the principal eigenvalue is in this case
\begin{align*}
	\lambda_1{\color{black}(\nu, \theta)} - \lambda_{\mathrm{exact}}{\color{black}(\nu, \theta)} = \frac{1}{540} I \nu  \left(2 \nu ^4-5 \nu ^3+5 \nu -2\right) \theta
   ^5+O\left(\theta ^6\right).
\end{align*}
\begin{figure}[t]
	\centering
	\includegraphics[scale=0.45]{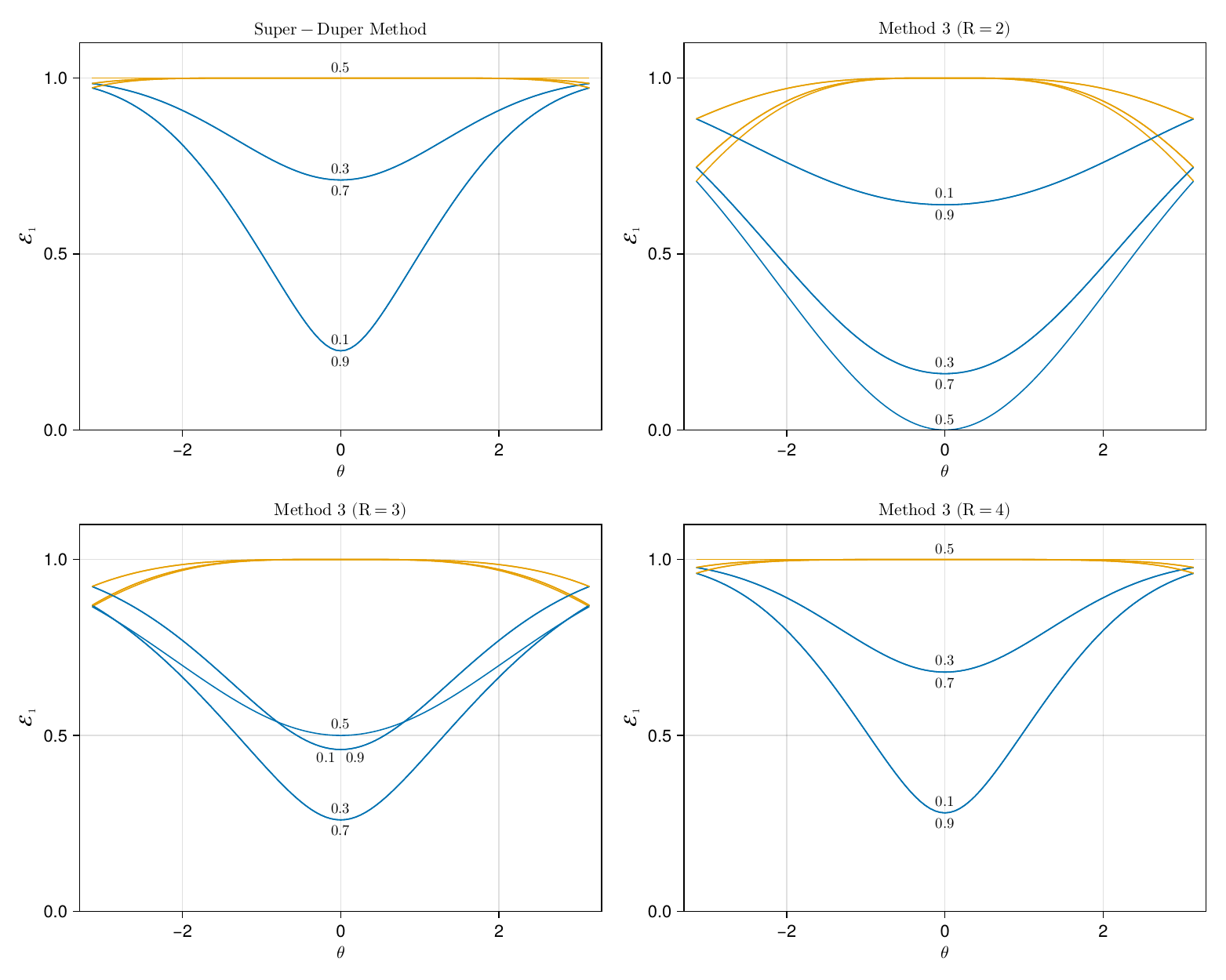}
	\caption{The dissipation error $\mathcal{E}_1$ with respect to $\theta$ for different values of $\nu$. The orange curves represent the error corresponding to the principal eigenvalue.}
	\label{fig:dissiaption}
\end{figure}
This error does not depend on any of the parameters. Hence, it is not possible to obtain a fifth-order correct eigenvalue with the right parameter choice. Nevertheless, if $\nu = 0$, $\nu = 1$, or $\nu = \frac{1}{2}$, then $\lambda_1(\nu, \theta) - \lambda_{\mathrm{exact}}(\nu, \theta) = 0$. For $\nu = 0$ and $\nu = 1$, this is expected, since the method is exact for these values of $\nu$, independently of the chosen parameters. If $\nu = \frac{1}{2}$, conditions \eqref{eq:cond_3} simplify to
\begin{align}
	S =  8-R, \quad T =  \frac{1}{8} (R-2)^2, \quad U =  \frac{1}{8} (R-6)^2.
	\label{eq:cond_ex}
\end{align}
It is a simple matter of calculation to show that for a scheme with parameters that satisfy \eqref{eq:cond_ex} and $\nu = \frac{1}{2}$, we have
\begin{align*}
	q_{i+\frac{1}{2}}^{n+2} = q_{i-\frac{1}{2}}^n \quad \text{and} \quad Q_i^{n+2} = Q_{i-1}^n.
\end{align*}
Therefore, conditions \eqref{eq:cond_3} yield a scheme that is exact for $\nu = \frac{1}{2}$. Alternatively, imposing conditions \eqref{eq:cond_2} with $S = 8-R$ also leads to a scheme that is exact for $\nu = \frac{1}{2}$, since these conditions reduce to conditions \eqref{eq:cond_ex}.
The principal eigenvalue of such schemes matches the exact eigenvalue for $\nu = \frac{1}{2}$, 
\begin{align}
	\lambda_1{\color{black}\left(\frac{1}{2}, \theta\right)} = \mathrm{e}^{-I \frac{\theta}{2}}.
    \label{eq:eigenvalue_12}
\end{align}
The spurious eigenvalue is
\begin{align*}
	\lambda_2{\color{black}\left(\frac{1}{2}, \theta\right)} = -\mathrm{e}^{-I \frac{\theta}{2}}.
\end{align*}
That these schemes are exact for $\nu = \frac{1}{2}$ is not a coincidence but a consequence of the amplification matrix. The elements of the amplification matrix consist of terms
\begin{align*}
	c\mathrm{e}^{-nI\theta}
\end{align*}
with $n \in \mathbb{N}_0$ and $c \in \R$. If the amplification matrix is in $\mathbb{R}$ (in this case we also call it the amplification factor), for example, for a finite volume method, then the eigenvalue of $A$, which in this case is just $A$, can only be of the form $\mathrm{e}^{-nI\theta}$ multiplied by some real factor. Hence, an eigenvalue like \eqref{eq:eigenvalue_12} is not possible. If $A \in \R^{2 \times 2}$, then the eigenvalues are the zeros of the characteristic polynomial, which is, in this case, a quadratic polynomial. Therefore, square root terms can appear in the eigenvalues, and these root terms are the reason why eigenvalues consisting of terms $\mathrm{e}^{-\frac{n}{2}I\theta}$ multiplied by some real factor can appear. Based on this argument, we expect it to be possible to modify higher-order schemes to obtain methods that are exact for different values of $\nu$.

\begin{figure}[t]
	\centering
	\includegraphics[scale=0.45]{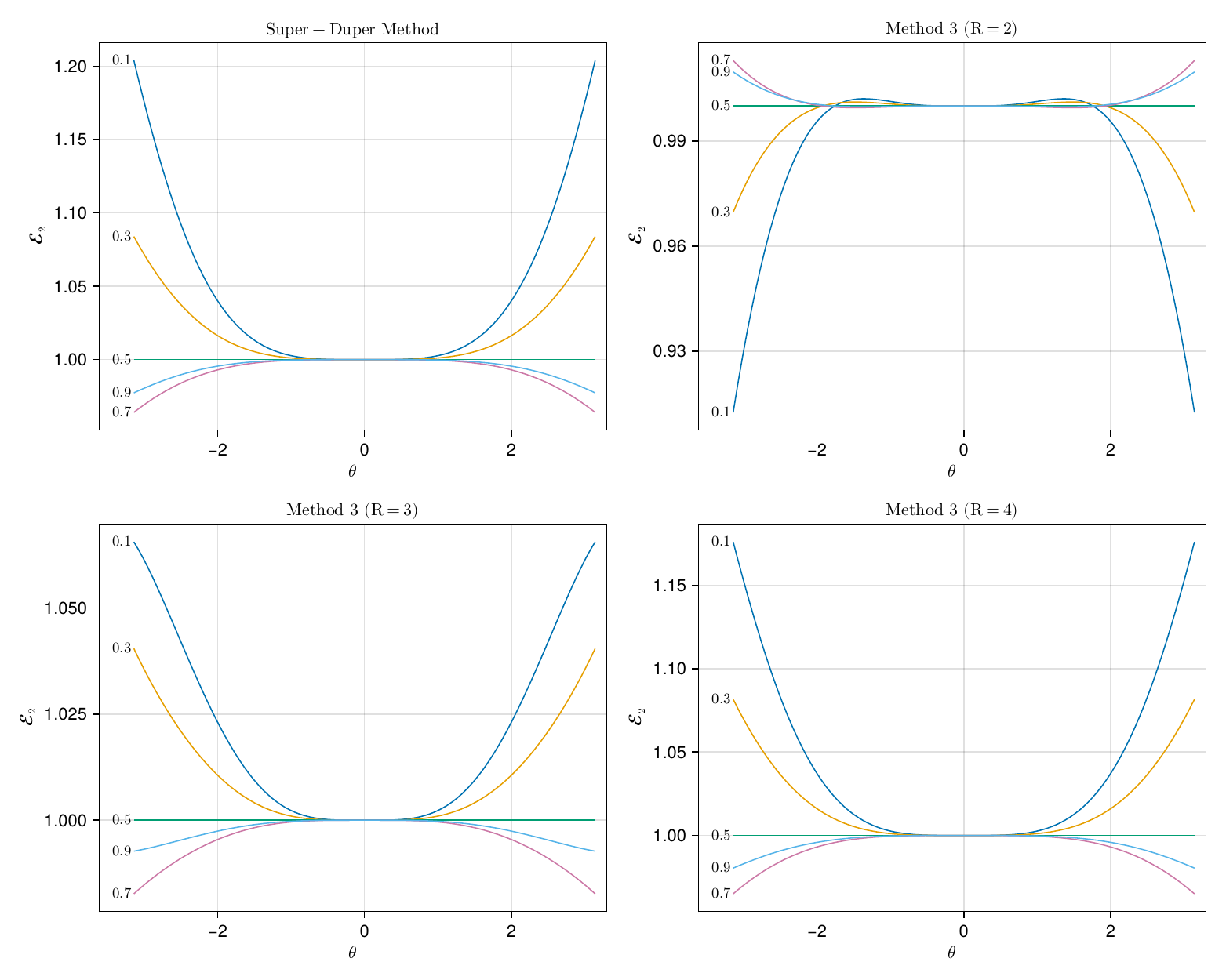}
	\caption{The first component of the dispersion error $\mathcal{E}_2$ corresponding to the principal eigenvalue with respect to $\theta$ for different values of $\nu$.}
	\label{fig:dispersion}
\end{figure}

Next, we analyze the dissipation and dispersion behavior of the schemes for certain parameter choices. Given the amplification matrix, we can define the dispersion and dissipation error in the spirit of \cite{Hirsch}. The dissipation error (or relative amplitude) $\mathcal{E}_{1}$ is defined as
	\begin{align*}
		\mathcal{E}_{1}(\nu, \theta) = \left(\frac{|\lambda_1(\nu, \theta)|}{|\lambda_\mathrm{exact}(\nu, \theta)|}, \frac{|\lambda_2(\nu, \theta)|}{|\lambda_\mathrm{exact}(\nu, \theta)|} \right)^\mathrm{T},
	\end{align*}
	and the dispersion error (or relative wave speed) $\mathcal{E}_2$ is defined as 
\begin{align*}
	\mathcal{E}_2(\nu, \theta) = \left(\frac{\phi_1(\nu, \theta) }{\phi_{\mathrm{exact}}(\nu, \theta)}, \frac{\phi_2(\nu, \theta) }{\phi_{\mathrm{exact}}(\nu, \theta)}\right)^\mathrm{T},
\end{align*} 
where $\phi_j(\nu, \theta)$ is the phase of the eigenvalue $\lambda_j(\nu, \theta)$ with  $j = 1,2$ and $\phi_{\mathrm{exact}}(\nu, \theta)$ is the phase of $\lambda_\mathrm{exact}(\nu, \theta)$.
We will focus our analysis on the schemes that have an at least third-order accurate principal eigenvalue, starting with schemes that have a fourth-order correct eigenvalue. Imposing conditions \eqref{eq:cond_3} results in eigenvalues 
\begin{align*}
	\lambda_{1,2}{\color{black}(\nu, \theta)} = & \ \frac{2 (2 \nu -1) \nu  e^{-I \theta }}{\nu +1}-4 \nu -\frac{6}{\nu -2}-2 \pm\frac{3 (\nu -1) \nu  e^{-2 I \theta }}{(\nu -2) (\nu +1) } \\
    & \ \sqrt{e^{3 I \theta } \left(3 I (2 \nu -1) \sin (\theta )+(1-2 \nu )^2 \cos
   (\theta )-4 (\nu -2) (\nu +1)\right)}.
\end{align*}
\begin{figure}[t]
	\centering
	\includegraphics[scale=0.45]{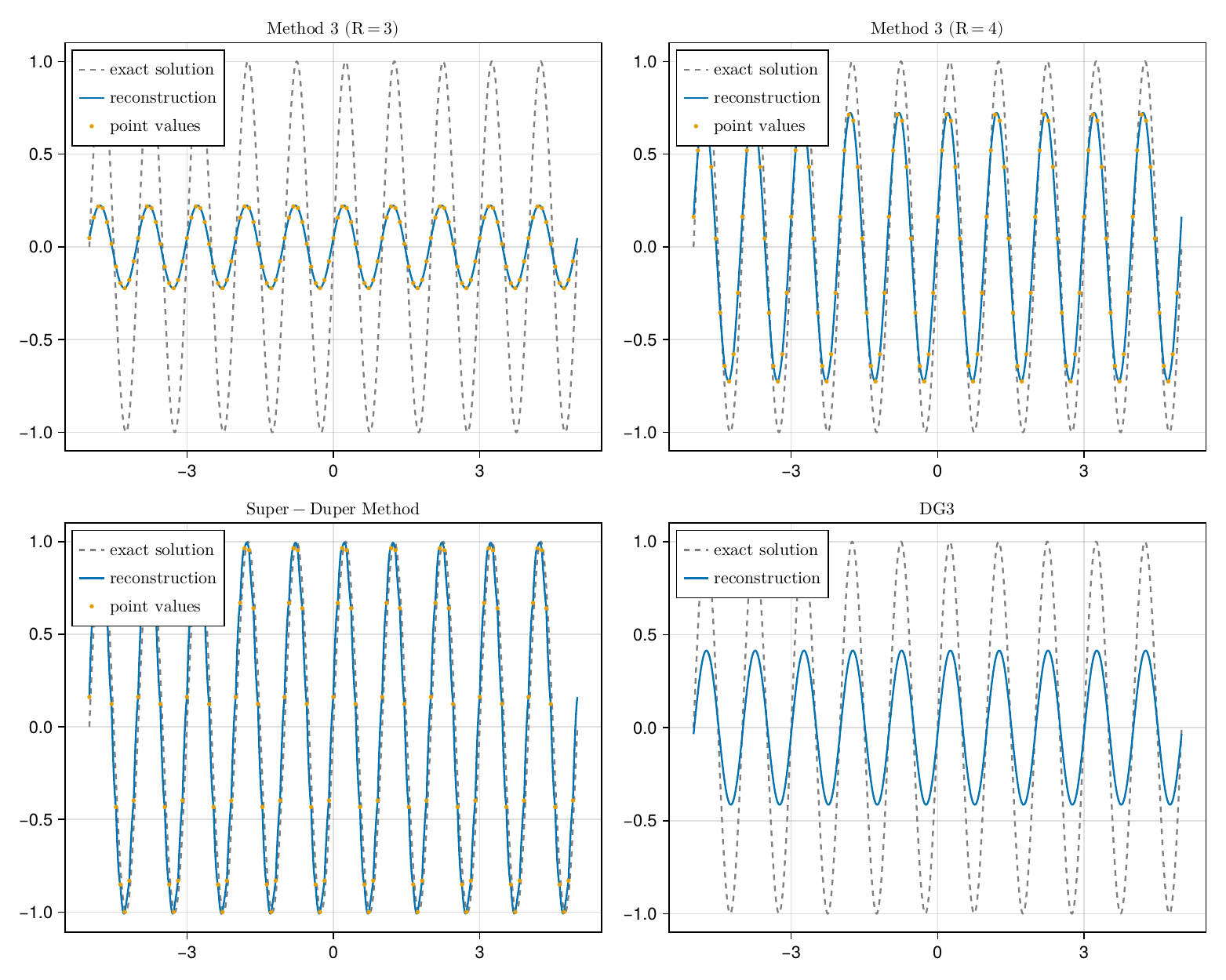}
	\caption{The numerical solutions computed with the different methods at $t = 2000$ in comparison to the analytical solution.}
	\label{fig:sine_wave}
\end{figure}

The dissipation and dispersion errors are always the same, no matter the choice of $R$. A useful way to determine $R$ is to make the expressions for the other parameters as simple as possible. The choice 
\begin{align*}
	R = \frac{6}{2 - \nu}
\end{align*}
accomplishes this, yielding
\begin{align*}
	S = \frac{6}{1 + \nu} \quad \text{and} \quad T = U = \frac{1}{2}.
\end{align*}
Phil Roe called this scheme the Super-Duper Method.
Imposing conditions \eqref{eq:cond_2} yields the eigenvalues
\begin{align*}
	\lambda_{1,2}\color{black}(\nu, \theta) =& \ 1+\frac{1}{6} \nu  \bigg(\nu ^2 (R+S)-e^{-I \theta } ((\nu -2) (\nu -1) (R+S)-6)-(R+S)-6\\
	& \ \pm \frac{1}{R+S}\bigg((\nu -1) e^{-2 I \theta } \bigg(e^{3 I \theta } (R+S)^3 (3 I (2 \nu -1) (R+S) \sin (\theta )\\
	& \ +\cos (\theta ) ((2 (\nu -1) \nu +5) (R+S)-36)\\
    &\ -2 (\nu -2) (\nu +1) (R+S) +36) \bigg)^{\frac{1}{2}} \bigg)\bigg).
\end{align*}
Both eigenvalues depend only on $(R+S)$. We set $R = S$ and consider the resulting schemes as a class depending on a single parameter $R$. We refer to the class of schemes with third-order accurate eigenvalues as Method 3. 

The dissipation error of the Super-Duper Method and Method 3 with $R = 2$, $R = 3$ (traditional Active Flux) and $R = 4$ with respect to $\theta$ can be seen in Figure \ref{fig:dissiaption} for different values of $\nu$. In  Figure \ref{fig:dispersion}, one can see the dispersion error of the principal eigenvalue with respect to $\theta$. A method with perfect dissipative and dispersive behavior would have a dissipation and dispersion error of $1$. The Super-Duper Method and Method 3 with $R=4$ have improved dissipative behavior compared to the traditional third-order scheme. This advantage appears to come with a cost, since these schemes are more dispersive than the classical version of the Active Flux scheme. Method 3 with $R=2$ is more dissipative than the traditional Active Flux method. The relative wave speed for the different values of $\nu$ is in some sense the opposite of that of the traditional method. For example, for $\nu = 0.1$, the dispersion error of the traditional method is larger than one, and for Method 3 with $R = 2$ it is smaller than one, provided $|\theta|$ is sufficiently large. 

The parameter $R$ in Method 3 allows us to adjust the dissipative and dispersive behavior of the method. Starting with $R = 3$ and increasing $R$ gives us a method that has an improved dissipative behavior compared to the classical Active Flux method and decreasing $R$ results in a method that has an inferior dissipative behavior.

\begin{figure}[t]
	\centering
	\includegraphics[scale=0.45]{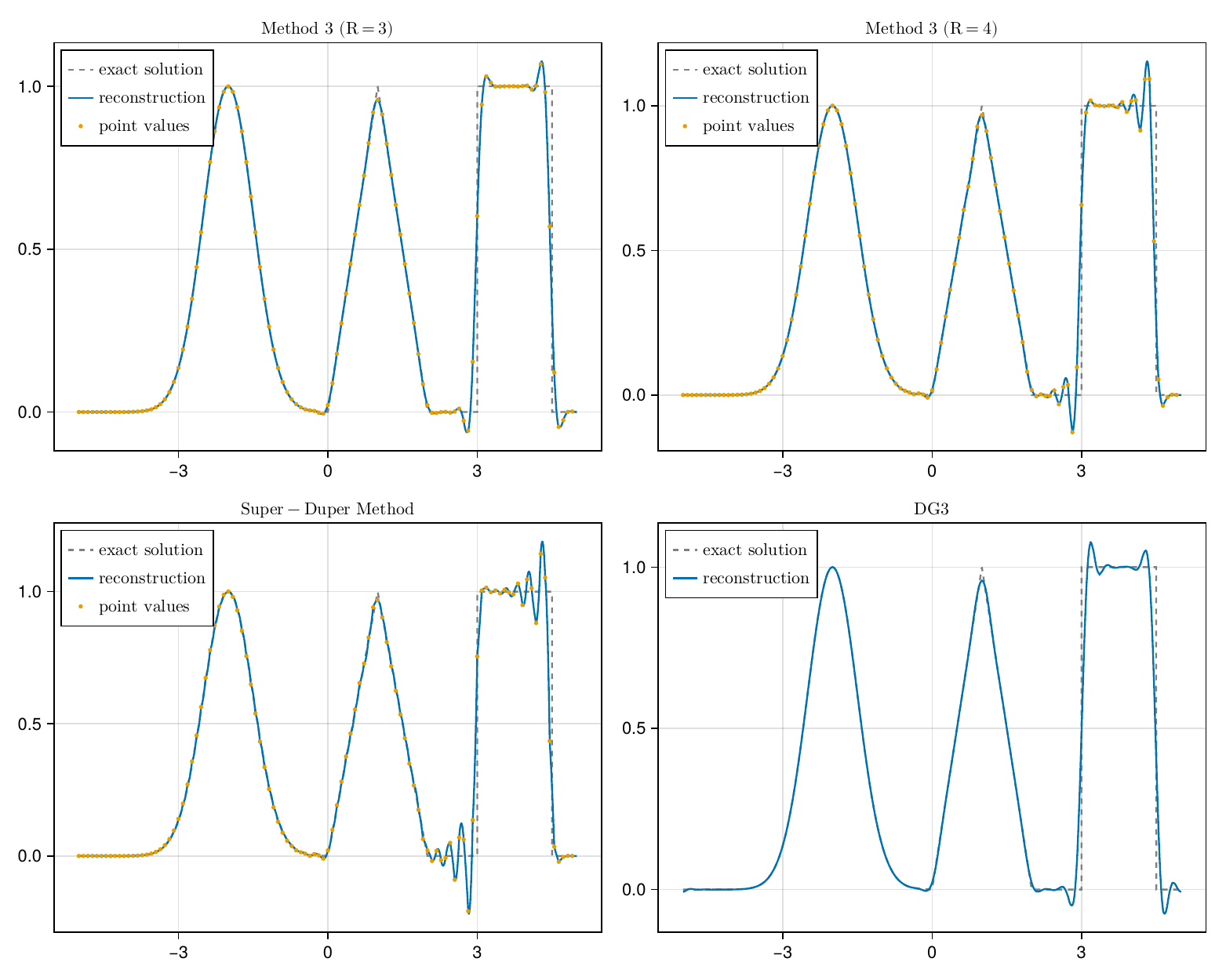}
	\caption{The numerical solutions computed with the different methods at $t = 20$ in comparison to the analytical solution.}
	\label{fig:shape_wave}
\end{figure}

\section{Numerical experiments}

For the numerical experiments, we choose $a = 1$, consider $[-5,5]$ as the computational domain divided into $110$ cells, and impose periodic boundary conditions.

For the first experiment, we use a high-frequency sine wave as initial data and compute numerical solutions with the Super-Duper Method and Method 3 with $R = 3$ and $R = 4$ at $t = 2000$. In addition, we compute a numerical solution with a third-order modal discontinuous Galerkin method (DG3). In our experiments, we found that the parameters do not impact the restrictions on the time step size. All Active Flux schemes considered in the experiments are stable for $\nu \leq 1$. To ensure that the comparison between the methods is as fair as possible, $\nu$ is chosen as approximately 95\% of the stability limit, i.e., $\nu = 0.2$ for DG3 and $\nu = 0.95$ for the Active Flux methods. The numerical solutions computed with the different methods are depicted in Figure \ref{fig:sine_wave}. The Super-Duper Method retains $100$ percent of the original amplitude, which explains its name.

\begin{figure}[t]
	\centering
	\includegraphics[scale=0.45]{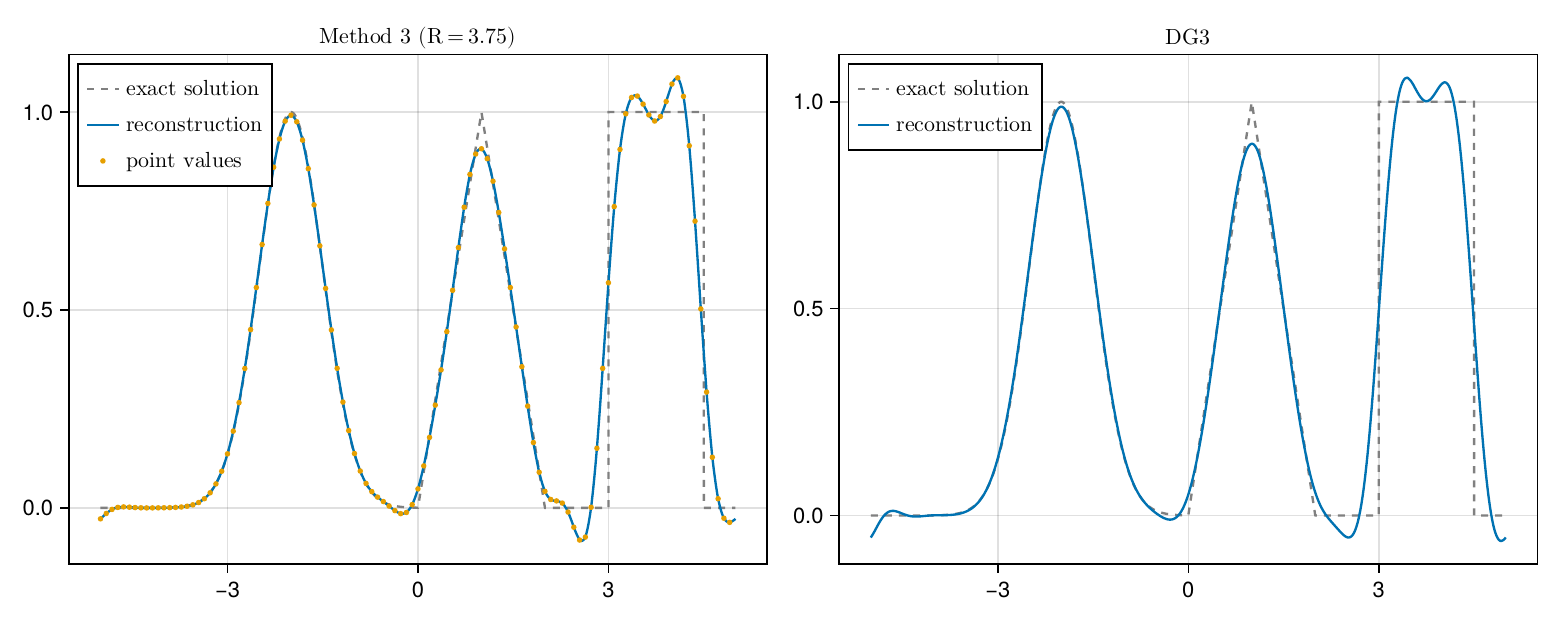}
	\caption{The numerical solutions computed with Method 3 ($R=3.75$) and DG3 at $t=1000$ in comparison to the analytical solution.}
	\label{fig:shape_wave_improved}
\end{figure}

Figure \ref{fig:shape_wave} shows the numerical results for a problem with initial data consisting of three different shapes. The methods with improved dissipative behavior and, as a consequence, worsened dispersive behavior produce oscillations near the discontinuities. Experimentally, values of $R \in [3,4]$ can be determined for which the oscillations remain acceptable. For example, in Figure \ref{fig:shape_wave_improved}, the square shape of the numerical solution computed with Method 3 ($R=3.75$) is approximately as accurate as that of the DG3 solution. In addition, the other two shapes of the numerical solution computed with Method 3 ($R=3.75$) are slightly more accurate than those of the DG3 solution.

\begin{figure}[t]
	\centering
	\includegraphics[scale=0.45]{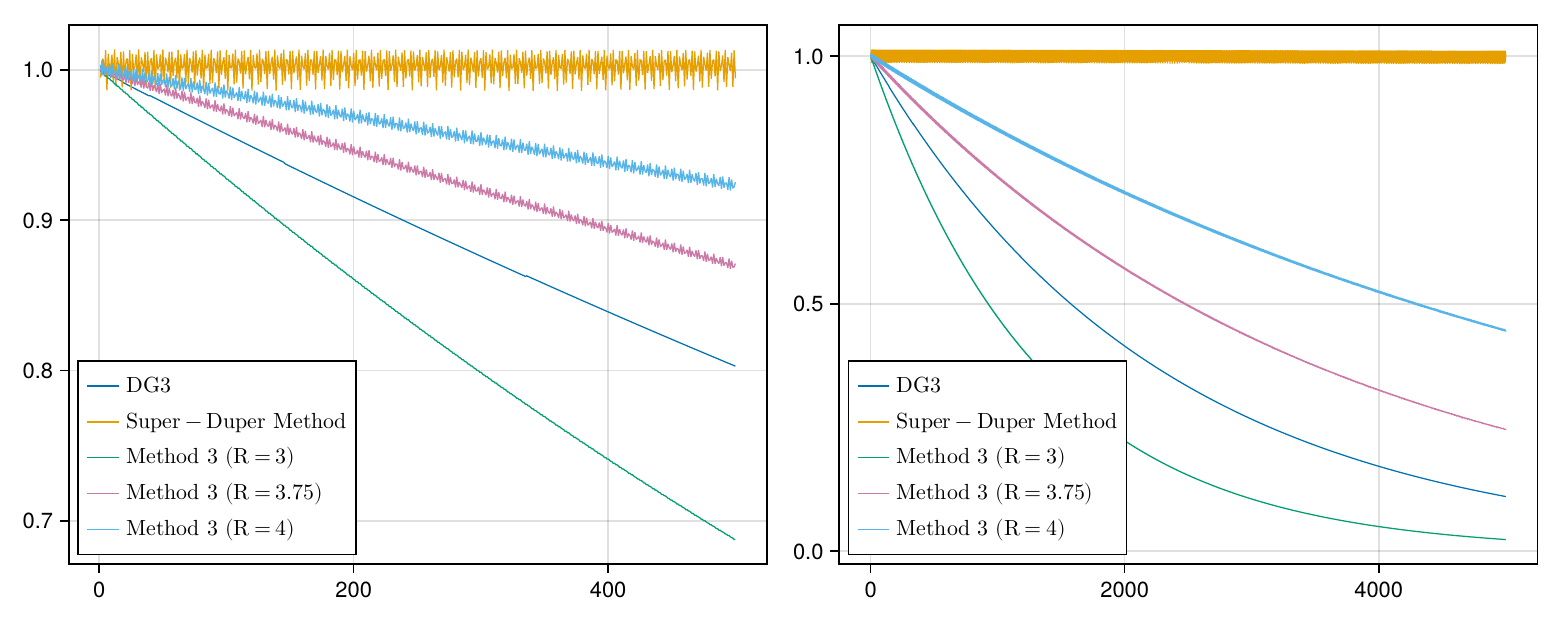}
	\caption{The maximum amplitude of the reconstruction computed with the different methods as a function of time $t$.}
	\label{fig:amplitude}
\end{figure}

\begin{figure}[t]
	\centering
	\includegraphics[scale=0.45]{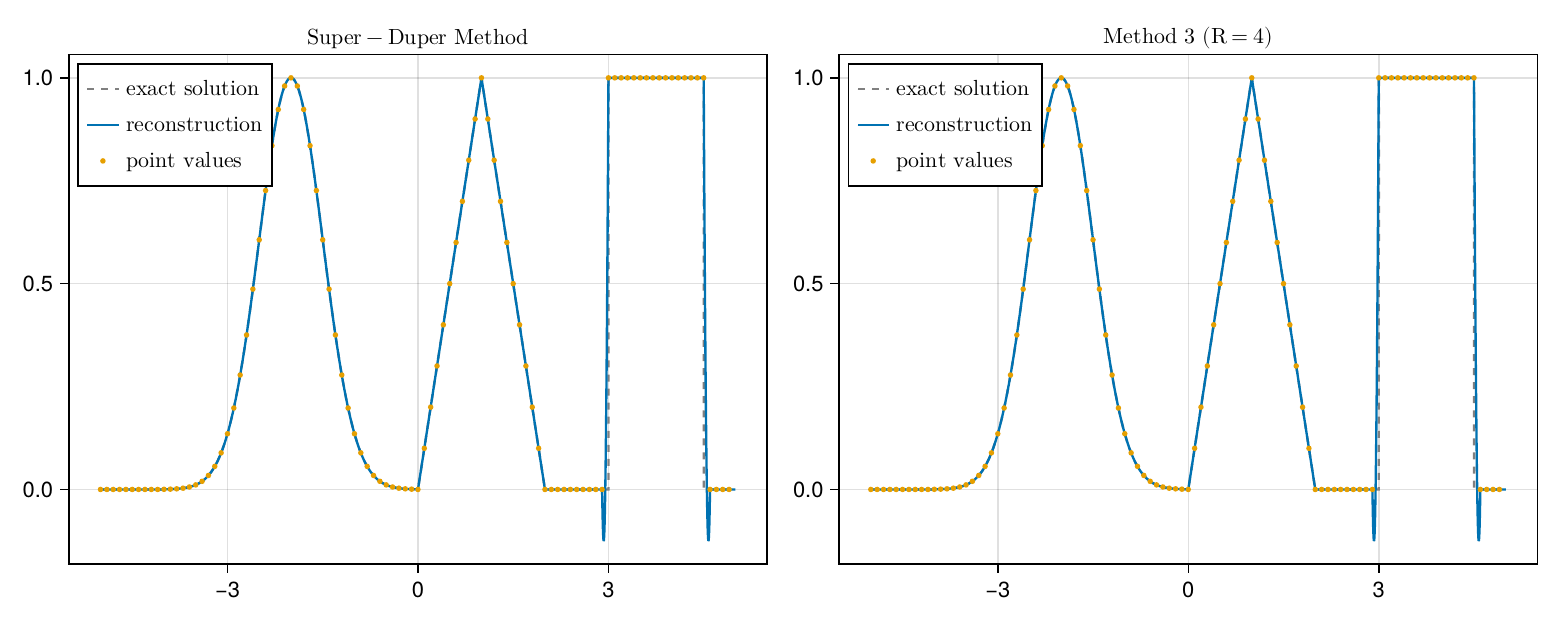}
	\caption{The numerical solution of the Super-Duper Method and Method 3 ($R=4$) with $\nu =0.5$ at $t = 10$ in comparison to the analytical solution.}
	\label{fig:shape_wave_2}
\end{figure}

Next, we again consider the high-frequency sine wave as initial data. In Figure \ref{fig:amplitude}, the maximum amplitude of the reconstruction obtained with the different methods is depicted. (The amplitude of the reconstruction is measured at every time step.) Even after $5000$ time steps, the amplitude of the numerical solution computed with the Super-Duper Method remains unchanged. Finally, as expected theoretically, the Super-Duper Method and Method 3 with $R = 4$ are exact for $\nu = 0.5$ (see Figure \ref{fig:shape_wave_2}).

\section{Conclusion and outlook}

In this work, we presented a reformulation of the traditional third-order Active Flux advection scheme as an upwind scheme with parameter-dependent correction terms. This formulation enables the selection of parameters, leading to a scheme with significantly improved dissipative behavior (the Super-Duper Method) and a class of schemes depending on a single parameter that allows control over the dissipative behavior (Method 3). These results were verified by numerical experiments. 

We have shown, using a simple problem, that the properties of a numerical method can be improved by introducing parameters and selecting them appropriately. 
In our numerical experiments this modification appeared to have no significant impact on the runtime or the computational complexity. 
A topic for future work is to extend the idea of optimizing numerical methods with the help of parameters to other problems, such as linear systems in one and two spatial dimensions and potentially to nonlinear conservation laws. Furthermore, it would be interesting to apply this idea to higher-order Active Flux methods.

\subsubsection*{Acknowledgement}
We acknowledge Lisa Lechner and Wasilij Barsukow for inspiring and productive discussions.

\end{document}